\documentstyle[12pt,amstex,intlim]{article}
\input tcilatex

\title{ON THE RIEMANN-HILBERT PROBLEMS}

\author{G. Giorgadze\\
Institute of Cybernetics\\
Georgian Academy of Sciences\\
e-mail: giorgadze@@rmi.acnet.ge}

\date{}

\begin{document}

\maketitle

\begin{abstract}
We discuss some topological aspects of the Riemann-Hilbert transmission problem
and Riemann-Hilbert monodromy problem on Riemann surfaces. In particular, we describe
the construction of a holomorphic vector bundle starting from the given representation
of the fundamental group and investigate the local behaviour of connexions on this
bundle. We give formul\ae\ for the partial indices of the Riemann-Hilbert transmission problem
in the three-dimensional case in terms of the correspoding vector bundle on the Riemann sphere.
\end{abstract}

\section*{Introduction}

The problem which we below call Riemann-Hilbert transmission problem (RHTP)
consists in finding piecewise holomorphic matrix functions which satisfy
certain transmission condition on the unit circle. It was first formulated
by B.~Riemann as an auxiliary proposition for solving the following problem:
starting with the given $m$ points on the Riemann sphere ${\bf CP}^1$ and $m$
nondegenerate matrices, one has to construct  a system of ordinary
differential equations (ODE's)  for which the given points are poles of the
first kind, and matrices are monodromy matrices around the given points.
Such systems are called Fuchsian and this problem is called the Riemann-Hilbert
monodromy problem (RHMP). Riemann himself did not give a detailed solution.
A valuable contribution to solving Riemann's problem was made by D.~Hilbert.
He solved the problem in certain particular case and stated an analogous
problem for the system of ODE's  which later became known as the Hilbert's 21-st
problem.

Plemelj in the paper ``Riemannsche Funktionenscharen mit gegebenen
Monodromiegruppe'' solved a weak version of Hilberts's 21-st problem. He
proved existence of a system of ODE's  which has the prescribed
monodromy, but also has a higher order pole. A system of this kind is called a
regular system. Plemelj used a different method than Hilbert. He made use of
a Cauchy type integral for solving the transmission problem. Besides that, he
brought a regular system to Fuchsian type by certain tranformation and
therefore gave a seemingly complete solution of RHMP.

However his proof was not completely rigorous, as was indicated by many
authors: [Mus], [V], [Kohn], [Ar-Il], [An-Bl]; this refers to RHTP as well as
to RHMP.

In subsequent years, research of the RHTP and RHMP was going indepedently.
RHTP was developed by the Georgian mathematical school and results of this
research were presented in the monographies of N.~Muskhelishvili [Mus] and
N.~Vekua [V]. Further results were obtained by W.~Koppelman [Kop], F.~Gachov
[Gah], Yu.~Rodin [Rod], G.~Khimshiashvili [Kh] and so on.

The history of the RHMP is even more interesing. It was Lappo-Da\-ni\-lev\-ski\u i
[LD] who was able to construct a Fuchsian system via its monodromy group.
Although Plemelj's proposition on the reduction of a regular system to
Fuchsian system was not true in general, Lappo-Da\-ni\-lev\-ski\u i's result remains
true because he considered the particular case when the monodromy matrices are
close to the identity matrix.

V.~Golubeva, A.~Bolibruch, V.~Leksin considered RHMP on $n$-dimensional
complex manifolds. They obtained many interesing results [Gol], [Lek1], [Lek2]
and their algebraic topological approach turned out to be useful for the
investigation of the monodromy representation of braid groups,
Knizhnik-Zamolodchikov and Yang-Baxter equations [Lek2], [Kohno]. Eventually,
A.~Bolibruch proved  that Plemelj's result is not true, i.~e. the answer on
the Hilbert's question on the existence of a Fuchsian system with prescribed
monodromy group is negative.

More detailed historical information may be found in the books
[Mus], [V], [Gah], [An-Bl]. In the sequel we will freely refer to these books.

The author became interested in these problems, because in last 20 years,
they found applications in mathematical physics, in particular, in holonomic
field theory [S-M-J], two dimensional Yang-Mills theory [At-Bot] and
topological quantum field theory [Wi].

The moduli space of holomorphic vector bundles on the Riemann surface
traditionally plays important role in complex analysis. Recently it also
became crucial for topological quantum field theory. It turned out that for
its description the language of loop groups is quite adequate, which
represents a generalization of ideas by Birkhoff and Bojarski to the case of
a Riemann surface of higher genus.

The holomorphic structure for a bundle on the Riemann sphere is determined
completely by the partial indices of the corresponding matrix function.
These are obtained via Birkhoff factorization of invertible matrix
functions. The condition of stability of a matrix function and the
corresponding notion of stability for vector bundles are expressible in the
language of partial indices developed by B.~Bojarski and I.~Gohberg with
M.~Krein. B.~Bojarski also explored the topology of stable matrix functions.
The diagonal matrix entering in the Birkhoff factorization of a matrix
function  represents a cocycle which defines a holomorpic vector bundle on $%
{\bf CP}^1.$ This bundle has connexions with regular singularities. The space
of Fuchsian connexions is contained in the space of connexions with regular
singularities. A.~Bolibruch discovered that for every conformal structure on
the marked Riemann sphere there exists a bundle which has no Fuchsian
connexion. He also explored the solution in a neighbourhood of a regular
singularity by Levelt's theory and transition from local solution to
global solutions. Due to this, the description of the moduli space of the
holomorphic vector bundles on a Riemann surface of higher genus became
possible in terms of local numeric invariants of ODE's. Our paper is devoted
to the same topic and written in the same spirit.

This text emerged from a talk given by the author at Prof. Bojarski's
seminar in the Institute of Mathematics of the Polish Academy of Sciences.
During discussions with Prof. B.~Bojarski several new ideas appeared which
will be considered in future publications.

I am very grateful to Prof. B.~Bojarski and the Banach Centre for the
invitation and warm hospitality.

\section{Riemann-Hilbert transmission problem}

Let $U^{+}$ be a bounded domain in the extended complex plane ${\bf CP}^1=%
{\bf C}^1\cup \infty ,$ with the boundary $L,$ $U^{-}$ -- the complementary
domain, so that $L=\overline{U}^{+}\bigcap \overline{U}^{-},{\bf CP}%
^1=U^{+}\bigcup L\bigcup U^{-}.$ For convenience of notation assume $0\in
U^{+}$ and $L$ is a piecewise smooth Jordan curve.

Let $\varphi (t)$ be a given function on $L,$ bounded everywhere on $L$ with
the possible exception of a finite number of points $s_1$, $s_2$,...,$s_m$, where
it satisfies the H\"older condition 
$$
\left| \varphi (t)\right| \leq \frac C{\left| t-s_j\right| ^\alpha }, 
$$
for any $s_j,j=1,2,...,m.$ The numbers $C,\alpha $ are positive constants and 
$\alpha <1.$

Consider the Cauchy type integral 
$$
\Phi (z)=\frac 1{2\pi i}\int_L\frac{\varphi (t)}{t-z}dt,\eqno{(1.1)} 
$$
where $z\in {\bf CP}^1$ and $t\in L$. The function $\Phi (z)$ is piecewise
holomorphic in ${\bf CP}^1\setminus L$ and for sufficiently large
$\left| z\right| $ we have a decomposition of $\Phi (z)$ into the sum
$$
\Phi (z)\sim \sum_{j=1}^\infty \frac{a_j}{z^j}\eqno{(1.2)} 
$$
where $a_j=-\frac 1{2\pi i}\int_Lt^{j-1}\varphi (t)dt,j=1,2,...$. Therefore,
we have $\Phi (\infty )=0.$

Let $t_0\in L$ and $t_0\notin \left\{ s_1,s_2,...,s_m\right\} ,$ then in
ordinary sense integral (1.1) does not exist. Consider the formal expression
$$
\Phi (t_0)=\frac 1{2\pi i}\int_L\frac{\varphi (t)}{t-t_0}dt.\eqno{(1.3)} 
$$

Describe about $t_0$ as centre a circle with so small a radius $\varepsilon $
that it intersects $L$ in two points $t_1$ and $t_2$. Denote by $l$ the
arc $t_1t_2.$ If for $\varepsilon \rightarrow 0$ the integral
$$
\frac 1{2\pi i}\int_{L\setminus l}\frac{\varphi (t)}{t-t_0}dt\eqno{(1.4)} 
$$
tends to a definite limit, then this limit is called the principal value of the
Cauchy type integral (1.1).

{\bf Remark. }It is obvious that, if the integral (1.3) exists in the
ordinary (i.~e. Riemann) sense, then the principal value also exists (but not
conversely). The integral (1.3) exists in the ordinary sense, if the
integral (1.4) tends to a definite limit whatever the arc $l$ cut off
around $t_0$ may be, as long as the length of this arc tends to zero; it is
essential for the definition of the principal value that the ends $t_1$ and
$t_2$ of the arc $l$ lie at equal distances from $t_0.$

{\bf Theorem 1.1. }If in the neighbourhood of $t_0$ the function $\varphi
(t) $ satisfies the H\"older continuity condition then there exists the
principal value of the Cauchy type integral.

This theorem allows to define the value of the function $\Phi (z)$ on the curve 
$L.$

Denote
$$
\Phi ^{\pm }(t_0)=\lim _{z\rightarrow t_0,z\in U^{\pm }}\Phi (z),t_0\in L, 
$$
then we have Sohocki\u i-Plemelj formulas for boundary values
$$
\Phi ^{+}(t_0)=\frac 12\varphi (t_0)+\frac 1{2\pi i}\int_L\frac{\varphi (t)}{%
t-t_0}dt, 
$$
$$
\Phi ^{-}(t_0)=-\frac 12\varphi (t_0)+\frac 1{2\pi i}\int_L\frac{\varphi (t)%
}{t-t_0}dt. 
$$

The aforementioned definitions and propositions also extend for those cases
when $\varphi $ is an $n$-vector or an $n\times n-$matrix function.

The Riemann-Hilbert transmission problem for vector valued functions can be
formulated as follows:

\ 

{\bf Riemann-Hilbert transmission problem:} Suppose one has H\"older
continuous matrix-functions $G:L\rightarrow GL(n,{\bf C})$. Describe the
totality of the piecewise holomorphic vector functions $\Phi (t)$ in $U^{+}\cup U^{-},$
which admit continuous boundary values on $L,$ satisfy the transmission condition
$$
\Phi ^{+}(t)=G(t)\Phi ^{-}(t),\eqno{(1.5)} 
$$
for any $t\in L$, and have finite order at $\infty .$

\ 

This problem is reduced to the system of singular integral equations in
the following manner.

Let the identity
$$
\Phi ^{+}(t)-\Phi ^{-}(t)=\varphi (t) 
$$
be satisfied on the curve $L.$ The piecewise holomorphic function, which has
finite order at $\infty $ is represented by the Cauchy type integral
$$
\Phi (z)=\frac 1{2\pi i}\int_L\frac{\varphi (t)}{t-z}dt+\gamma (z),%
\eqno{(1.6)} 
$$
where $\gamma (z)=(\gamma _1(z),\gamma _2(z),...,\gamma _n(z))$ is the
principal part of the function $\Phi (z)$ at $\infty .$ Therefore, our
problem is reduced to finding functions $\varphi $ and $\gamma .$ Let
us use Sohocki\u i-Plemelj formulas for $\Phi (z)$ and substitute the obtained
expression in the transmission condition (1.5). We obtain the system of
singular integral equations for the functions $\varphi _1,\varphi
_2,...,\varphi _n:$
$$
A(t_0)\varphi (t_0)+\frac{B(t_0)}{\pi i}\int_L\frac{\varphi (t)}{t-t_0}%
dt=F(t_0),\eqno{(1.7)} 
$$
where $A(t_0)={\bf 1+}G(t_0),B(t_0)={\bf 1-}G(t_0),F(t_0)=(G(t_0)-{\bf 1}%
)\gamma (t)$ and ${\bf 1}$ is the unit matrix in $GL(n,{\bf C})$.

The function $F(t_0)$ contains an unknown polynomial $\gamma (t_0).$ It must be
chosen so that the system (1.7) must have solutions. This last condition
is satisfied if and only if the system
$$
\int_Lf(t)\psi ^k(t)dt=0,k=1,2,...,k^{^{\prime }}\eqno{(1.8)} 
$$
is satisfied, where $\psi ^k(t),k=1,2,...,k^{^{\prime }}$ forms a complete
system of linearly independent solutions of the adjoint of problem (1.5).

{\bf Properties of the solutions of RHTP:}

{\bf Property 1. }If $\Phi _1(z),...,\Phi _n(z)$ are solutions of the RHTP,
then for any polynomial functions
$p_1(z),...,p_n(z)$ the function 
$$
\Phi_1(z)p_1(z)+...+\Phi_n(z)p_n(z) 
$$
again is a solution.

{\bf Property 2.} Suppose for some $z_0$ we have $\Phi (z_0)=0.$ Then
$$
\frac{\Phi (z)}{z-z_0} 
$$
again is a solution of the RHTP.

{\bf Property 3.} Let $k$ be the number of linearly indepedent solutions of
the system of singular integral equations (1.7), then the order at infinity
of nontrivial solutions of the RHTP is at least $k$ (note, that $%
k\geq$ind$G+(m+1)n-k^{^{\prime }},$ where $m$ is order of the pole at
infinity).

Let
$$
\Psi (z)=\left( 
\begin{array}{cccc}
\Phi _1^1(z) & \Phi _1^2(z) & ... & \Phi _1^n(z) \\ 
\Phi _2^1(z) & \Phi _2^2(z) & ... & \Phi _2^n(z) \\ 
... & ... & ... & ... \\ 
\Phi _n^1(z) & \Phi _n^2(z) & ... & \Phi _n^n(z) 
\end{array}
\right) 
$$
where the column $\Psi ^i(z)=(\Phi^i_1(z),...,\Phi^i_n(z))$ of the matrix
$\Psi (z)$ consists of solutions of the RHTP and let $\det \Psi (z)$ be not
identically zero. In this case the system of the solutions
$\Psi ^1(z),...,\Psi ^n(z)$ is called a fundamental system of solutions.

If $\det \Psi (z)\neq 0$ for any $z\in {\bf C}$ (including the curve $L$),
then $\Psi^1(z)$,...,$\Psi^n(z)$ is called normal. Note, that
$\det \Psi(t)\neq 0$ for $t\in L$ means, that
$\det \Psi^+(t)\neq 0,\det \Psi^-(t)\neq 0.$

Determinant of the system of normal solutions may be 0 at infinity, or
$\infty $ may be a pole for it. Denote by $\kappa _1,\kappa _2,...,\kappa _n$
the order at infinity of the system of normal solutions $\Psi
^1(z),...,\Psi ^n(z)$ and consider the matrix-function $\Psi (z)z^{d_K},$
where
$$
d_K=\left( 
\begin{array}{ccc}
z^{\kappa _1} &  & 0 \\  
& ... &  \\ 
0 &  & z^{\kappa _n} 
\end{array}
\right)\eqno{(1.9)} 
$$
is a diagonal matrix with entries $z^{\kappa _i}$ and the integers
$\kappa_1,\kappa _2,...,\kappa _n$ satisfies the conditions
$$
\kappa _1\geq \kappa _2\geq ...\geq \kappa _n.\eqno{(1.10)} 
$$

It is clear that $\lim _{z\rightarrow \infty }\det (\Psi
(z)z^{d_K})=c<\infty $ and therefore $\Psi (z)z^{d_K}$ is a holomorphic
matrix function in the neighbourhood of $\infty .$

The normal system of solutions is called canonical, if the function
$\Delta (z)=\det (\Psi (z)z^{d_K})$ is not zero at $\infty$. Denote by
$\chi (z)$ the system of canonical solutions.

The integers (1.10) are called partial indices of the homogeneous
transmission problem or of the matrix function $G(t).$ The number
$$
\kappa =\frac 1{2\pi i}\Delta _L\arg \det G(t) 
$$
is called the global index or simply index of the RHTP.

{\bf Theorem 1.2.}(Muskhelishvili [Mus]). For every $G(t)$ canonical
solutions always exist. The sequence (1.9) does not depend on the
considered canonical solution and
$$
\kappa =\kappa _1+\kappa _2+...+\kappa _n. 
$$

The matrix-function
$$
\chi _0(z)=\Psi (z)z^{d_K} 
$$
is holomorphic invertible in $U^{-}$ and $\det $ $\chi _0(z)\neq 0,$ for
every points $t\in L.$ Using the identity (1.5) we obtain the
Birkhoff factorization of the matrix-function $G(t):$
$$
G(t)=\chi ^{+}(t)t^{d_K}\left[ \chi _0^{-}(t)\right] ^{-1}. 
$$

It is reasonable to detalize this Birkhoff theorem as we'll use it again.

Denote by $\Omega $ the space of all H\"older-continuous loops $%
G:L\rightarrow GL(n,{\bf C}).$ It is a Banach Lie group with natural
norm and operetion.

Let

$\Omega ^{+}=\{f\in \Omega :f$ is the boundary value of the matrix function
holomorphic in U$^{+}\},$

$\Omega ^{-}=\{f\in \Omega :f$ is the boundary value of the matrix function
holomorphic in U$^{-}$ and is regular at infinity $f(\infty )={\bf 1\}.}$

{\bf Theorem 1.3. }(Birkhoff) Any loop $f\in \Omega $ can be represented as
$$
f(t)=f^{-}(t)d_Kf^{+}(t),\eqno{(1.11)} 
$$
where $f^{\pm }\in \Omega ^{\pm }$ and $d_K$ is a diagonal loop (1.9) with
condition (1.10).

The diagonal matrix $d_K$ will be called the characteristic loop of the
corresponding matrix-function, whereas $K=(k_1,k_2,...,k_n)$ -- the
characteristic multiindex. Two loops $f,g\in \Omega $ will be called
equivalent, if $f$ and $g$ have identical characteristic multiindices. For
$K=(k_1,k_2,...,k_n)$, denote by $\Omega_K$ the set of equivalence classes
of loops $\Omega $ and call it the Birkhoff-Bojarski stratum. The topological
structure of $\Omega_K$ has been studied by Bojarski, who showed, that if
$L$ is a Jordan curve, then the strata $\Omega_K$ are connected. The
representation (1.11) is not unique, but if one fixes $f^{+}$ (or $f^{-})$
then $f^{-}$ (respectively $f^{+})$ will be uniquely defined.

The Banach Lie group $\Omega ^{+}\times \Omega ^{-}$ acts analytically on $%
\Omega $ via
$$
f\stackrel{\alpha }{\longmapsto }h_1fh_2^{-1},f\in \Omega ,h_1\in \Omega
^{+},h_2\in \Omega ^{-}. 
$$
It is clear, that the orbit of the diagonal matrix $d_K$ by the action $%
\alpha $ is $\Omega _K.$

{\bf Theorem 1.4 }(Disney [Dis]) The stability subgroup $H_K$ of $%
f $ under the action $\alpha $ consists of those pairs $(h_1,h_2)$ of upper
triangular matrix-functions where the $(i,j)$-th entry in $h_1$ is a
polynomial in $z$ of degree at most $(k_1-k_2)$ and
$f=h_1fh_2^{-1}$, the space $H_K$ has finite dimension
$$
\dim H_K=\stackunder{k_i\geq k_j}{\sum }(k_i-k_j+1). 
$$

The stratum $\Omega _K$ is a locally closed analytical submanifold
of $\Omega $ and codimension of $\Omega _K$ in $\Omega $ is equal to
$$
\dim \Omega /\Omega _K=\stackunder{k_i\rangle k_j}{\sum }(k_i-k_j-1). 
$$

Consider the holomorphic vector bundle on ${\bf CP}^1$ which is obtained
by the covering of the Riemann sphere ${\bf CP}^1$ by three open sets
$\{U^{+},U^{-},U_3={\bf CP}^1\backslash \{0,\infty\}\}$, with transition
functions
$$
g_{13}=h_1:U^{+}\cap U_3\rightarrow GL(n,{\bf C}), 
$$
$$
g_{23}=h_2d_K:U^{-}\cap U_3\rightarrow GL(n,{\bf C}). 
$$
It is denoted by $E\rightarrow {\bf CP}^1.$ From the Birkhoff theorem it
follows, that every holomorphic vector bundle splits into direct sum of
the line bundles
$$
E\cong E(k_1)\oplus ...\oplus E(k_n).\eqno{(1.12)} 
$$

{\bf Remark.} Possibility of decomposition of a holomorphic vector bundle into
the sum (1.12) is proved by A. Grothendieck, without applying the Birkhoff
theorem.

The numbers $k_1$,...,$k_n$ are the Chern numbers of the line bundles
$E(k_1)$,...,$E(k_n)$ and satisfy the conditions $k_1\geq$...$\geq k_n.$ The
integer-valued vector $K=(k_1,...,k_n)\in {\bf Z}^n$ is called the splitting
type of the holomorphic vector bundle $E.$ It defines uniquely the
holomorphic type of the bundle $E.$

Connection between partial indices $\kappa _1,...,\kappa _n$ of the RHTP,
characteristic multiindex of the matrix-function $f\in \Omega $ and splitting
type of the holomorphic vector bundle $E$ are presented in the following
summarizing theorem:

{\bf Theorem 1.5.} There is a one-to-one correspodence between the strata $%
\Omega _K$ and holomorphic vector bundles on ${\bf CP}^1.$

Denote by $O(E)$ the sheaf of germs of holomorphic sections of the bundle $E,$
then the solutions of the RHTP are elements of the zeroth cohomology group
$H^0({\bf CP}^1,O(E)),$ therefore the number $l$ of the linearly independent
solutions is $\dim H^0({\bf CP}^1,O(E)),$ as the Chern number $c_1(E)$ of the
bundle $E$ is equal to index $\det G(t),$ we obtained the known criterion of
solvability of the RHTP. In particular the following theorem is true:

{\bf Theorem 1.6.} The RHTP has solutions if and only if $c_1(E)\geq 0$ and
the number $l$ of linearly independent solutions is
$$
l=\dim H^0({\bf CP}^1,O(E))=\sum_{i=1}^nk_i+1. 
$$

Suppose $L$ is the same as above and $G:L\rightarrow GL(n,{\bf C})$ is a
discontinuous function  of the first kind at the point $s_1\in L,$ i.~e.
$\stackunder{t\rightarrow s_1+0}{\lim }G(t)\neq \stackunder{t\rightarrow s_s-0}{\lim }%
G(t).$ Denote by $G(s_1+0)=\stackunder{t\rightarrow s_1+0}{\lim }%
G(t),G(s_1-0)=\stackunder{t\rightarrow s_1-0}{\lim }G(t).$

Let
$$
G=G^{-1}(s_1+0)G(s_1-0)\eqno{(1.13)} 
$$
and $\Gamma =\frac 1{2\pi i}\ln G,$ so that if $\lambda _i$ are eigenvalues of 
$G,$ then $\mu_i=\frac 1{2\pi i}\ln \lambda _i$ satisfies the conditions $%
0\leq Re\mu _i<1.$

Consider the functons
$$
\omega ^{+}(z)=(z-s_1)^\Gamma \stackrel{def}{=}e^{\Gamma \ln (z-s_1)}, 
$$
$$
\omega ^{-}(z)=\left( \frac{z-s_1}{z-z_0}\right) ^\Gamma \stackrel{def}{=}%
e^{\Gamma \ln \left( \frac{z-s_1}{z-z_0}\right) }, 
$$
where $z_0$ is some fixed point in $U^{+}.$ It is known that $\omega ^{+}(z)$
is a single valued matrix-function on ${\bf C}\setminus l_1,$ where $l_1$ is a
curve with endpoints $s_1$ and $\infty .$ $\omega ^{-}(z)$ is a single valued
matrix-function on ${\bf C}\setminus l_2,$ where $l_1$ is a curve with
endpoints $z_0$ and $s_1.$

Suppose
$$
\stackunder{t\rightarrow s_1+0}{\lim }(t-z_0)^\Gamma =(s_1-z_0)^\Gamma , 
$$
then
$$
\stackunder{t\rightarrow s_1-0}{\lim }(t-z_0)^\Gamma =e^{2\pi
i}(t_1-z_0)^\Gamma =G(t_1-z_0)^\Gamma.\eqno{(1.14)} 
$$

Introduce new vector-functions $f_1^{+}(z)$ and $f_1^{-}(z):$
$$
f_1^{+}(z)=(z-s_1)^\Gamma G^{-1}(s_1+0)f^{+}(z), 
$$
$$
f_1^{-}(z)=\left( \frac{z-s_1}{z-z_0}\right) ^{-\Gamma }f^{-}(z). 
$$
They are holomorphic respectively on $U^{\pm }$ and satisfy the
transmission condition
$$
G(s_1+0)(z-s_1)^\Gamma f^{+}(z)=G(t)\left( \frac{z-s_1}{z-z_0}\right)
^\Gamma f^{-}(z), 
$$
or
$$
f_1^{+}(z)=(z-s_1)^{-\Gamma }G^{-1}(s_1+0)G(t)\left( \frac{z-s_1}{z-z_0}%
\right) ^\Gamma f_1^{-}(z). 
$$

Let us denote $G_1(z)=(z-s_1)^{-\Gamma }G^{-1}(s_1+0)G(t)\left( \frac{z-s_1}{%
z-z_0}\right) ^{-\Gamma }$ and prove that $G_1(t)$ is continuous at the point 
$s_1.$ Indeed
$$
G_1(s_1+0)=\stackunder{t\rightarrow s_1+0}{\lim }\left[ (t-s_1)^{-\Gamma
}G^{-1}(s_1+0)G(t)\left( \frac{t-s_1}{t-z_0}\right) ^\Gamma \right]
=(s_1-z_0)^{-\Gamma }. 
$$
To calculate $G_1(s_1-0)$ we use (1.13), (1.14) and obtain

$G_1(s_1-0)=\stackunder{t\rightarrow s_1-0}{\lim }\left[ (t-s_1)^{-\Gamma
}G^{-1}(s_1+0)G(t)\left( \frac{t-s_1}{t-z_0}\right) ^\Gamma \right] =$

$=\stackunder{t\rightarrow s_1-0}{\lim }\left[ (t-s_1)^{-\Gamma
}GG^{-1}\left( \frac{t-s_1}{t-z_0}\right) ^\Gamma \right]
=(s_1-z_0)^{-\Gamma }.$

Now consider the general case. Let $s_1,...,s_m\in L$ be points of discontinuity
and let there be finite limits $G(s_j+0)=\stackunder{%
t\rightarrow s_1+0}{\lim }G(t)$ and $G(s_j-0)=\stackunder{t\rightarrow s_1-0%
}{\lim }G(t).$ The curve $L$ will be understood to be the union of smooth,
nonintersecting arcs $L_1,L_2,...,L_m$ with definite positive directions.
Therefore ends of arcs $L_j$ (j=1,2,...,m) are $s_j$ and $s_{j+1}.$

Suppose

$G^j=G^{-1}(s_j+0)G(s_j-0)$ and $\Gamma ^j=\frac 1{2\pi i}\ln G^j,$ so that
if $\lambda _j^i$ are eigenvalues of $G^j,$ then $\mu _i^j=\frac 1{2\pi i}\ln
\lambda _j^i.$ Denote $\rho _i^j=Re\mu _i^j$ and normalize the choice of $%
\ln $ demanding that $0\leq \rho _i^j<1.$

Consider the matrix-functions
$$
\Omega _j^{+}(z)=A_jG(s_j+0)(z-s_j)^{\Gamma _j},\Omega _j^{-}(z)=B_j\left( 
\frac{z-s_j}{z-z_0}\right) ^{\Gamma _j}, 
$$
where $A_j,B_j$ are constant matrices:
$$
A_1=E,A_j=\left[ \prod_{k=1}^{j-1}\Omega _k^{+}(s_j)\right] ^{-1}, 
$$
$$
B_1=E,B_j=\left[ \prod_{k=1}^{j-1}\Omega _k^{-}(s_j)\right] ^{-1},
j=2,3,...m. 
$$

The functions $\Omega _j^{+}(z)$ are holomorphic respectively in $U^{\pm
}.$

Introduce new vector-functions
$$
f^{+}(z)=\prod_{j=1}^m\Omega _j^{+}(z)f_1^{+}(z),
$$
$$
f^{-}(z)=\prod_{j=1}^m\Omega _j^{-}(z)f_1^{-}(z).
$$
Use the transmission condition (1.5) and obtain:
$$
f_1^{+}(t)=\left[ \left( \prod_{j=1}^m\Omega _j^{+}(t)\right)
^{-1}G(t)\prod_{j=1}^m\Omega _j^{-}(t)\right] f_1^{-}(t). 
$$

{\bf Proposition 1.1.}The matrix-function
$$
G_1(t)=\left( \prod_{j=1}^m\Omega _j^{+}(t)\right)
^{-1}G(t)\prod_{j=1}^m\Omega _j^{-}(t) 
$$
is continuous at points $s_1,...,s_m.$

{\bf Proof.} Calculate $G_1(s_j+0)$ and $G_1(s_j-0).$
$$
G_1(s_j+0)=
$$
$$
=\stackunder{t\rightarrow s_1+0}\lim\left(
\prod_{k=j+1}^m\Omega _k^{+}(s_j)\right) ^{-1}(t-s_j)^{-\Gamma
_j}G^{-1}(s_j+0)A_j^{-1}\left( \prod_{k=1}^{j-1}\Omega _k^{+}(s_j)\right)
^{-1}=
$$
$$
=G(s_j+0)\prod_{k=1}^{j-1}\Omega _k^{-}(s_j)B_j\left( \frac{t-s_j}{t-z_0}%
\right) ^{\Gamma _j}\prod_{k=j+1}^m\Omega _k^{-}(s_j)= 
$$
$$
=\left( \prod_{k=j+1}^m\Omega _k^{+}(s_j)\right) ^{-1}\stackunder{%
t\rightarrow s_1+0}{\lim }\left( (t-s_j)^{-\Gamma _j}\left( \frac{t-s_j}{%
t-z_0}\right) ^{\Gamma _j}\right) \prod_{k=j+1}^m\Omega _k^{-}(s_j)= 
$$
$$
=\left( \prod_{k=j+1}^m\Omega _k^{+}(s_j)\right) ^{-1}(s_j-z_0)^{-\Gamma
_j}\prod_{k=j+1}^m\Omega _k^{-}(s_j), 
$$
$$
G_1(s_j-0)=\stackunder{t\rightarrow s_1-0}{\lim }\left(
\prod_{k=j+1}^m\Omega _k^{+}(s_j)\right) ^{-1}(t-s_j)^{-\Gamma
_j}G^{-1}(s_j+0)A_j^{-1}\left( \prod_{k=1}^{j-1}\Omega _k^{+}(s_j)\right)
^{-1} 
$$
$$
G(s_j-0)\prod_{k=1}^{j-1}\Omega _k^{-}(s_j)B_j\left( \frac{t-s_j}{t-z_0}%
\right) ^{\Gamma _j}\prod_{k=j+1}^m\Omega _k^{-}(s_j)= 
$$
$$
=\left( \prod_{k=j+1}^m\Omega _k^{+}(s_j)\right) ^{-1}\stackunder{%
t\rightarrow s_1-0}{\lim }\left( (t-s_j)^{-\Gamma _j}G_jG_j^{-1}\left( \frac{%
t-s_j}{s_j-z_0}\right) ^{\Gamma _j}\right) \prod_{k=j+1}^m\Omega
_k^{-}(s_j)= 
$$
$$
=\left( \prod_{k=j+1}^m\Omega _k^{+}(s_j)\right) ^{-1}(s_j-z_0)^{-\Gamma
_j}\prod_{k=j+1}^m\Omega_k^-(s_j). 
$$

Therefore, the RHTP with discontinuity points reduces to the transmission
problem considered at the very beginning, but in this case it is necessary
to find solutions, which are holomorphic respectively in $U^{\pm }$ and its
boundary values have discontinuity points. It can be proved that if
$(f$ $^{+},f^{-})$ is a solution of RHTP with discontinuity points, then this
solution extends continuously to these points too. It means, that there is a
system of canonical solutions $\chi_0(z)$ of the transmission problem, which
satisfies the following conditions:

1. $\det \chi (z)\neq 0,$ on ${\bf C}$ with possible exception of points $%
s_1,s_2,...,s_m.$

2. There are the diagonal matrix-function $d_K,$ that then
$$
\stackunder{z\rightarrow \infty }{\lim }\chi (z)d_K(z)=c\neq 0<\infty . 
$$

3. If $s_j$ is some singular point, then
$$
\stackunder{z\rightarrow s_j}{\lim }(z-s_j)^\varepsilon \chi (z)=0, 
$$
for some real number $\varepsilon >0$.

\section{The system of ordinary differential equations with regular
singularity}

Consider the system on a small disk $U\subset ${\bf C} with center 0, 
$$
\frac{df}{dz}=A(z)f(z),f(z)=(f^1(z),...,f^n(z))\in {\bf C}^n,\eqno{(2.1)} 
$$
where $A(z)$ is a holomorphic matrix function on $U^{*}=U\backslash \{0\}$.

Let $p:\widetilde{U^*}\to U^*$ be the universal covering of $U^*$ and let $%
\xi$ and $z$ denote the local coordinate on $\widetilde{U^*}$ and $U^*$,
respectively.

The system (2.1) has $n$ linearly independent holomorphic solutions in a small
neighbourhood of $z_0\in U^*$. Denote the space of solutions  by $\Re$. If $%
f\in \Re$, then $f$ is a holomorphic function on $\widetilde{U^*}$. Let $%
\Gamma$ be the group of deck transformations of the covering 
$$
p:\widetilde{U^*}\to U^*. 
$$

If $\alpha\in\Gamma$, then $\alpha$ defines the automorphism $%
\alpha^*:\Re\to\Re$ of the solution space in this manner: 
$$
\alpha^*f=f\circ\alpha^{-1},\ {\rm i. e.}\
(\alpha^*f)(\xi)=f(\alpha^{-1}\xi). 
$$

Clearly, $\alpha ^{*}f$ is also a solution to (2.1) and therefore a map 
$$
\rho :\Gamma \to GL(n,{\bf C}),\ \alpha \longmapsto \alpha ^{*}.\eqno{(2.2)} 
$$
is obtained.

If $\beta\in\Gamma$ is another element, then $(\alpha\beta
)^*=\alpha^*\beta^*$, i.~e. the map (2.2) is a homomorphism. Thus 
$$
f=(f\circ\alpha)\rho(\alpha).\eqno{(2.3)}
$$

The homomorphism $\rho$ is called the monodromy representation corresponding
to the system (2.1).

Let $\Phi(z)$ be the fundamental system of solutions to (2.1) and let $\Phi
_1(z)$ be another invertible solution of the matrix ODE's: 
$$
\frac{d\Phi _1}{dz}=A(z)\Phi_1(z). 
$$
Then $\Phi _1(z)=\Phi (z)G$ with some constant matrix $G\in GL(n,{\bf C)}$.
Instead of (2.3) we get $\Phi _1(z)=(\Phi _1\circ \alpha )\rho _1(\alpha )$
with some 
$$
\rho _1:\Gamma \to GL(n,{\bf C}). 
$$
So 
$$
\Phi (z)G=(\Phi (z)G\circ \alpha )\rho _1(\alpha )=(\Phi \circ \alpha )G\rho
_1(\alpha ). 
$$

But $\Phi (z)=(\Phi (z)\circ \alpha )\rho (\alpha ),$ thus $(\Phi (z)\circ
\alpha )\rho (\alpha )G=(\Phi (z)\circ \alpha )G\rho _1(\alpha ).$ Hence $%
\rho _1(\alpha )=G^{-1}\rho (\alpha )G,$ where $G$ is the same for all $%
\alpha .$ We see that to a system (2.1) there corresponds a class of
mutually conjugate representations $\rho :\Gamma \to GL(n,{\bf C)}$. We will
call this class the monodromy representation or simply monodromy.

The group of deck transformations $\Gamma $ is now the infinite cyclic group
generated by the deck transformation $\alpha $ which corresponds to one trip
around $0$ counterclockwise. Clearly, $\ln \xi $ is a holomorphic function
on $\tilde U^{*}$ and $\ln (\alpha \xi )=\ln \xi +2\pi i.$ Let $G=\rho
(\alpha ^{-1})$ so that 
$$
\Phi (\alpha \xi )=\Phi (\xi )G.\eqno{(2.4)} 
$$

Let $E=\frac 1{2\pi i}\ln G,$ so that if $\lambda _j$ are eigenvalues of $G$
and $\mu _j$ of $E,$ then $\mu _j=\frac 1{2\pi i}\ln \lambda _j.$ Denote $%
\tau _j=Re$ $\mu _j$ and normalize the choice of $\ln $ demanding that $%
0\le\tau_j<1$.

Introduce the function $\xi ^E=e^{E\ln \xi }$ (which is holomorphic on $%
\tilde U^{*}$): 
$$
(\alpha \xi )^E=e^{E(\ln \xi +2\pi i)}=\xi ^EG. 
$$
Then by (2.4) 
$$
\Phi (\alpha \xi )(\alpha \xi )^{-E}=\Phi (\xi )GG^{-1}\xi ^{-E}=\Phi (\xi
)\xi ^{-E}. 
$$
Hence $\Phi (\xi )\xi ^{-E}$ can be considered as a single-valued
holomorphic function on $U^{*}.$

Consider any sector $\Sigma $ having vertex at $0$. 0 is a regular
singularity of this system, if for the covering sector $\Sigma ^{^{\prime }}$
on $U^{^{\prime }} $ and for any solution $f(\xi ),$ the restriction $f(\xi
)\mid_{\Sigma'}$ has at most polynomial growth as $z\to0$
remaining in $\Sigma .$

Analogically one defines the regular singularity of the $n$-th order
differential equation 
$$
x^{\left( n\right) }(z)+a_1(z)x^{\left( n-1\right)
}(z)+...+a_{n-1}(z)x^{^{\prime }}(z)+a_n(z)x(z)=0.\eqno{(2.5)} 
$$

Observe that the system (2.1) and the equation (2.5) have regular singular
points.

{\bf Theorem 2.1.(Poincar\`e) }Let $f(z)$ be some solution of the system
(2.1). Then $f(z)$ can be represented as follows:
$$
f(\xi )=Z(z)\xi ^E, 
$$
where $Z$ is holomorphic on $U^*$.

\ {\bf Proposition 2.1.} Every coordinate function $f_j(\xi )$ of a solution 
$f(\xi )$ is 
$$
f_j(\xi )=\sum_{p,q}\xi ^{\tau _p}h_{p,q}(z)\ln ^{l_q}\xi ,\eqno{(2.6)} 
$$
$$
0\leq Re\tau _q<1,l_q\in Z,l_q\geq 0. 
$$

Let $n_{pq}^j$ denote order of the zero of the function $h_{pq}(z)$ at
the point 0, and let $n^j=\min_{p,q}n_{p,q}^j.$ There is defined a map 
$$
\varphi:\Re\to{\bf Z},\ \varphi(f)=\min_{j=1,...,n}n^j. 
$$

The map $\varphi $ will be called Levelt's normalization of the solution $f(z).$
It has the following properties:

1. $\varphi(\lambda f)=\varphi(f)$, if $\lambda\in C^*$;

2. $\varphi(0)=\infty$;

3. $\varphi(f_1+f_2)\geq \min (\varphi (f_1),\varphi(f_2))$, with equality
if $\varphi (f_1)\neq \varphi(f_2)$.

From the algebraic viewpoint $\varphi $ is a nonarchimedean valuation on 
$\Re $ over the trivial valuation on ${\bf C}$.

The integer valued function $\varphi $ defines a filtration of $\Re :$ 
$$
0\subset \Re ^0\subset \Re ^1\subset \Re ^2\subset ...\subset \Re ^m\subset
\Re ,\eqno{(2.7)} 
$$
such that $\varphi $ is constant on the quotient space $\Re ^j/\Re ^{j-1}$
and if $k_j^i=\varphi (\Re ^j/\Re ^{j-1}),$ then $k_1>k_2>...>k_m.$ Let $%
d_j=\dim (\Re ^j/\Re ^{j-1})$. We say that $\varphi $ takes the value $k_j$
with multiplicity $d_j$.

We shall use also the notation 
$$
\varphi ^1=\varphi ^2=...=\varphi ^{d_1}=k_1>\varphi ^{d_1+1}=...=\varphi
^{d_1+d_2}=k_2>...>\varphi ^{d_1+d_2+...+d_{m-1}}%
$$
$$
=\varphi ^{d_1+d_2+...+d_m}=k_m.%
$$

Note that 
$$
\varphi^1\geq\varphi^2\geq...\geq\varphi^m. 
$$

By definition of $\varphi $ it follows that 
$$
\varphi (\alpha ^{*}f)=\varphi (f), 
$$
hence it follows that $\alpha^*$ preserves the filtration (2.7) and the
monodromy matrix $G$ is upper triangular.

A basis $f_1(z)$, $f_2(z)$,...,$f_n(z)$ of the solution space $\Re,$ 
satisfying the conditions $\varphi(f_i)=\varphi^i$ and such that the monodromy
matrix $G$ is upper triangular, will be called a Levelt's basis.

\ {\bf Theorem 2.2. [Le].} The fundamental system of solutions $\Phi (\xi )$
related to a Levelt's basis is 
$$
\Phi (\xi )=U(z)z^\Psi \xi ^E, 
$$
where $U(z)$ is holomorphic on $U^{*}$ and $\det U(z)\ne 0$, $\Psi =$ diag$%
(\varphi ^1,...\varphi ^m)$ and $E=\frac 1{2\pi i}G$ is upper triangular.

{\bf Remark.} If $0$ is a regular singular point, then $U(z)$ is meromorphic
in $U$, i.~e. $U(z)$ is a single valued function.

{\bf Proposition 2.2. }Let $U(z)$ be holomorphically invertible at $z=0$ and
let
$$
L(z)=\Psi +z^\Psi Ez^{-\Psi } 
$$
be holomorphic. Then the system
$$
df=\omega f 
$$
is Fuchsian at $z=0,$ where $\omega =\frac{d\Phi (z)}{dz}\Phi ^{-1}(z).$

{\bf Proof.}

$$
\frac{d\Phi (z)}{dz}=\frac{dU(z)}{dz}z^\Psi \xi ^E+\frac 1zU(z)z^\Psi E\xi
^E= 
$$
$$
\frac 1z(z\frac{dU(z)}{dz}+U(z)L(z))z^\Psi \xi ^E, 
$$
where
$$
L(z)=\Psi +z^\Psi Ez^{-\Psi}, 
$$
then
$$
\frac{d\Phi (z)}{dz}\Phi ^{-1}(z)=\frac 1z(z\frac{dU(z)}{dz}%
+U(z)L(z))U^{-1}(z) 
$$
is Fuchsian at point $0.$

If $0$ is a pole of order one for the matrix valued function $A(z),$ then
the system (2.1) is called Fuchsian. Let 
$$
A={\rm Res}_{z=0}A(z), 
$$
then (2.1) gives 
$$
\frac{df}{dz}=\frac Azf(z) 
$$

\ 

{\bf Proposition 2.3.} 1) Every Fuchsian system is regular.

2) $0$ is a regular singular point for the equation (2.5) if and only if the
functions $z^ja_j(z)$ are holomorphic at $0.$

\ {\bf Remark.} 1) The set of regular singular systems contains the set of
Fuchsian systems.

2) The ordinary differential equation (2.5) is regular if and only if it is
Fuchsian.

From the proposition 2.3 follows that the coefficients $%
a_1(z),a_2(z),...,a_n(z)$ are holomorphic in some punctured neighdourhood of 
$0$ and $a_1(z)$ has there at most a pole of the 1-st order, ..., $a_i(z)$
--- at most a pole of the $i$-th order, ..., $a_n(z)$ --- at most a pole of
the $n-$th order.

It turns out that (2.5) is regular at the point $0$ if and only if the
system is describing the behavior of the vector 
$$
(f^1(z),...,f^n(z))=(x(z),\frac{dx(z)}{dz},...,\frac{d^{n-1}x(z)}{dz^{n-1}}) 
$$
i.~e. the system 
$$
\frac{df(z)}{dz}=\left( 
\begin{array}{cccc}
0 & -1 & ... & 0 \\ 
... & ... & ... & ... \\ 
0 & 0 & ... & -1 \\ 
-a_n(z) & -a_{n-1}(z) & ... & -a_1(z) 
\end{array}
\right) f(z) 
$$
is regular at $0$. It is well known that (2.5) is regular at $z=0$ if and
only if it is Fuchsian at $z=0$.

The systems $\frac{df(z)}{dz}=A(z)f(z)$ and $\frac{dg(z)}{dz}=B(z)g(z)$ are
holomorphically (meromorphically) equivalent if there exists in a
neighbourhood of $0$ a holomorphic (meromorphic) at $0$ matrix function $%
H:V\to GL(n,{\bf C})$, such that the transformation $(z,f(z))=(z,H(z)f(z))$
maps one equation to another, i.~e. 
$$
B(z)=\frac{dH(z)}{dz}H^{-1}(z)+H(z)A(z)H^{-1}(z).\eqno{(2.8)} 
$$

If two systems of equations are equivalent then their monodromy groups are
conjugate. Besides, if $0$ is a regular singular point for the system, this
system is equivalent to $\frac{df}{dz}=\frac Azf(z),$ where $A$ is a constant
matrix.

To eliminate ambiguity we introduce some standard definitions.

Meromorphic connexion at the point $0$ is called a pair $(F,\nabla)$, where $%
F$ is an $n$-dimensional vector space over the field $K={\cal O}[\frac1z]$,
whereas $\nabla:F\to F$ is an operator which satisfies the Leibniz rule 
$$
\nabla(h,s)=\frac{dh}{dz}s+h\nabla s,%
$$
for each function $f\in K$ and $s\in F$.

Let $e_1,e_2,...,e_n$ be a basis of $F$ and let $\nabla e_i$ be expressed in
this basis in the following manner: 
$$
\nabla e_i=-\sum_{j=1}^n\theta _{ij}(z)e_j%
$$
where $\theta =(\theta _{ij}(z))\in$ End$(n,K)$, then for $%
s=\sum_{j=1}^ns_j(z)e_i$ we shall obtain 
$$
\nabla s=\sum_{i=1}^n(\frac{ds_i(z)}{dz}-\sum_{j=1}^n\theta
_{ij}(z)e_j(z))e_i.%
$$

By the last formula it follows, that $\nabla s=0$ is equivalent to $\frac{ds%
}{dz}=\theta s$, or $(d-\theta )s=0$.

Let us denote the matrix-valued 1-form $\theta dz$ by $\omega $, then the
system will be $(d-\omega )=0$ and the connexion will be $\nabla =d-\omega $%
. Connexions are gauge equivalent if and only if correspoding systems of
equations are equivalent, i.~e. satisfy the equality (2.8).

\section{Connection between RHTP and RHMP}

Let $G(t)$ be the transmission matrix-function for the RHTP and suppose it
is piecewise constant, i.~e. $G(t)=G_iG_{i-1}...G_1$ if $t\in s_is_{i+1}.$
Let $\chi_0(z)$ be the canonical solution and $\omega =d\chi_0(z)\chi_0^{-1}(z)$
be a form which is single-valued on the Riemann sphere and holomorphic
outside the points $s_1,...,s_m.$

{\bf Proposition 3.1.} The monodromy matrices of the regular system of
ODE's
$$
df=\omega f\eqno{(3.1)} 
$$
are $G_1,G_2,...,G_m.$

Indeed, consider RHTP with the piecewise constant transmission function $G(t)$
and denote by $\chi_0(z)$ the canonical solution of the corresponding
transmission problem. Let $z_0\in U^-.$ Take some singular point $s_i$
and let $\gamma_i$ be a loop beginning at $z_0$ and going around the singular
point $s_i$ along a small circle. It in transmission condition follows, that
the extension of $\chi_0^-(z)$ along the loop $\gamma _i$ goes to
$G_i\chi_0^-(z).$

{\bf Theorem 3.1. }(Plemelj) If some monodromy matrix is diagonalizable,
then the system (3.1) is Fuchsian.

{\bf Proof.} Let $Y(\xi )$ be the fundamental system of solutions and let the
monodromy matrix $E_j$ which corresponds to the singular point $s_j$ be
diagonal. By theorem 2.1 in the neighbourhood of $s_j$ one can represent
$Y_j(\xi )$ as follows:
$$
Y_j(\xi )=U_j(z)(\xi -s_j)^{E_j} 
$$
where $U_j(z)$ is holomorphic on $U^{*}.$ By Sauvage's lemma there exists a
matrix $\Gamma (z)$ holomorphically invertible outside $s_j$ such that
$$
\Gamma (z)U_j(z)=V_j(z)(z-s_j)^{\Psi _j} 
$$
where $V_j(z)$ is holomorphically invertible at $s_{j}$ and $\Psi
=diag(\varphi ^1,...,\varphi ^n).$ Introduce a new dependent variable
$g=\Gamma(z)f.$ By (2.8) we have%
$$
\frac{dg(z)}{dz}=\left( \frac{d\Gamma (z)}{dz}\Gamma ^{-1}(z)+\Gamma
(z)A(z)\Gamma ^{-1}(z)\right) g(z). 
$$

Therefore the new system is still Fuchsian outside $s_j.$

We want to prove that if conditions of the theorem are satisfied, then
the system is Fuchsian at the point $s_j$ too.

Take Levelt's $Y(\xi )$ (theorem 2.2) in the neighbourhood $s_j:$
$$
Y(\xi )=U(z)(z-s_j)^{\Psi _j}(\xi -s_j)^{E_j}. 
$$
We repeat the calculation presented during the proof of proposition 2.2.
$$
\frac{d\Phi (z)}{dz}=\frac{dU(z)}{dz}(z-s_j)^{\Psi _j}(\xi -s_j)^{E_j}+\frac
1zU(z)\Psi _j(z-s_j)^{\Psi _j}(\xi -s_j)^{E_j}+ 
$$
$$
+\frac 1zU(z)(z-s_j)^{\Psi _j}E_j(\xi -s_j)^{E_j}= 
$$
$$
=\frac 1z(z\frac{dU(z)}{dz}+V(z)L(z))(z-s_j)^{\Psi _j}(\xi -s_j)^{E_j}, 
$$
where
$$
L(z)=\bar \Psi +(z-s_j)^{\Psi _j}E_j(z-s_j)^{-\Psi _j}, 
$$
then we obtain:
$$
\frac{d\Phi (z)}{dz}\Phi ^{-1}(z)=\frac{1.}z(z\frac{dU(z)}{dz}%
+U(z)L(z))U^{-1}(z), 
$$
$E_j$ is diagonal, because $L(z)=\Psi _j+E_j$ is holomorphic. Therefore, by
proposition 2.3 our system is Fuchsian at the point $s_j$ too.

\section{Extension of a bundle with connexion}

Let $X$ be a Riemann surface of genus $g$ and $S=\left\{
s_1,s_2,...,s_m\right\} $ be a set of marked points on $X.$ Denote by $%
X_m=X\setminus S.$ Let $\tilde X\to X_m$ be the universal covering map of $%
X_m$, then it is a bundle with fibre $\pi _1(X_m,z_0),$ where $z_0\in X_m$. $%
\pi _1(X_m,z_0)$ is isomorphic to the group of deck transformations of this
covering and therefore acts on $\tilde X.$

Let 
$$
\rho:\pi _1(X_m,z_0)\to GL( n,{\bf C})\eqno{(4.1)}%
$$
be some representation.

Consider the trivial principal bundle $\tilde X\times GL(n,{\bf C})\to
\tilde X$ (or vector bundle $\tilde X\times {\bf C}^n\to \tilde X$). The
quotient space $\tilde X\times GL(n,{\bf C})/\sim$ gives the locally
trivial bundle on $X_m,$ where $\sim$ is an equivalence relation
identifying the pairs $(\tilde x,g)$ and $(\sigma \tilde x,\rho (\sigma )g),$
for every $\tilde x\in \tilde X,g\in GL(n,{\bf C})$ (or $g\in {\bf C}^n$).
Denote the obtained bundle by ${\bf P}_\rho \to X_m$ (or ${\bf E}_\rho \to
X_m$) and call it the bundle associated with the representation $\rho.$ In
obvious form this bundle according to the transformation functions may be
constructed in the following manner.

Let $\left\{U_\alpha\right\} $ be a simple covering of $X_m,$ i.~e. every
intersection $U_{\alpha _1}\cap U_{\alpha _2}\cap ...\cap U_{\alpha _k}$ is
connected and simply connected. For each $U_\alpha,$ we choose a point $%
z_\alpha \in U_\alpha $ and join $z_0$ and $z_\alpha $ by a $\gamma _\alpha $
starting at $z_0$ and ending at $z_\alpha$. For a point $z\in U_\alpha \cap
U_\beta $ we choose a path $\tau _\alpha \subset U_\alpha $ which starts at $%
z_\alpha $ and ends at $z.$ Consider 
$$
g_{\alpha \beta }\left( z\right) =\rho \left( \gamma _\alpha \tau _\alpha
\left( z\right) \tau _{\beta ^{}}^{-1}\left( z\right) \gamma _\beta
^{-1}\right) .\eqno(4.2) 
$$
We see that 
$$
g_{\alpha \gamma }\left( z\right) =g_{\beta \alpha }\left( z\right) 
$$
and 
$$
g_{\alpha \beta }g_{\beta \gamma }\left( z\right) =g_{\alpha \gamma }\left(
z\right) 
$$
on $U_\alpha \cap U_\beta \cap U_\gamma .$

The cocycle $\left\{ g_{\alpha \beta }\left( z\right) \right\}$ does not
depend on the choice of $z.$ Hence from this cocycle we obtain a flat vector
(or principal) bundle, which is denoted by ${\bf E}_\rho ^{\prime }$ $({\bf P%
}_\rho ^{^{\prime }}).$

Let $\left\{ t_\alpha \left( z\right) \right\} $ be a trivialization of our
bundle, i.~e. 
$$
t_\alpha :p^{-1}\left( U_\alpha \right) \to GL(n,{\bf C}) 
$$
is a holomorphic mapping. Consider the matrix valued $1$-form $\left\{ \omega
_\alpha \right\} :$ 
$$
\omega _\alpha =-t_\alpha ^{-1}dt_\alpha . 
$$
$\left\{ g_{\alpha \beta }\left( z\right) \right\} $ are constant on the
intersection $U_\alpha \cap U_\beta $ and $g_{\alpha \beta }(z)t_\beta
\left( z\right) =t_\alpha \left( z\right)$, so on $U_\alpha \cap U_\beta $
the identity $\omega _\alpha =\omega _\beta$ holds. Indeed, replacing $%
t_\beta $ by $t_\beta ^{-1}g_{\alpha \beta }$ in the expression $\omega
_\beta =-t_\beta ^{-1}dt_\beta$, we obtain 
$$
\omega _\beta =-t_\alpha ^{-1}g_{\alpha \beta }\left( z\right) dt_\alpha
g_{\alpha \beta }^{-1}\left( z\right) =-t_\alpha ^{-1}dt_\alpha .%
$$
So, $\omega =\left\{ \omega _\alpha \right\} $ is a holomorphic 1-form on $%
X_m$ and therefore is a connexion form of the bundle ${\bf P}_\rho ^{\prime
}\to X_m.$ Corresponding connexion is denoted by $\nabla ^{\prime }.$ We
will extend the pair $\left( {\bf P}_\rho ^{\prime },\nabla ^{\prime
}\right) $ to $X.$ As the required construction is of local character, we
shall extend ${\bf P}_\rho ^{\prime }\to X_m$ to the bundle ${\bf P}_\rho
^{\prime \prime }\to X_m\cup \left\{ s_i\right\} ,$ where $s_i\in S.$

First consider the extension of the principal bundle ${\bf P}_\rho
^{^{\prime }}\rightarrow X_m.$

Let a neighbourhood $V_i$ of the point $s_i$ meet $U_{\alpha _1},U_{\alpha
_2},...U_{\alpha _k}$. As we noted when constructing the bundle from
transition functions (4.1) only one of them is different from identity. Let
us denote it by $g_{1k}$, then $g_{1k}=G_i,$ where $G_i$ is the monodromy which
corresponds to the singular point $s_i.$ Mark a branch of the many valued
function $\left( \tilde z-s_i\right) ^{E_i}$ containing the point $\tilde
s_i\in \tilde U_i$ (where $E_i=\frac 1{2\pi i}\ln G_i$). Thus the marked
branch defines a function 
$$
g_{01}=(z-s_i)^{E_i}.\eqno{(4.3)} 
$$
Denote by $g_{02}$ the extension of $g_{01}$ along the path which goes around
$s_i$ counterclockwise, and similarly for other points. At last on
$U_i\cap U_{\alpha _k}\cap U_{\alpha _1}$ we shall have:
$$
g_{0k}(z)=g_{01}(z)G_i=g_{01}(z)g_{0k}(z).%
$$

The function $g_{0k}:V_i\to GL(n,{\bf C})$ is the one defined at the point
$s_i,$ and takes there value coinciding with the monodromy matrix. It means,
that we made extension of the bundle to the point $s_i.$ In a neighbourhood
of $s_i$ one will have
$$
\omega _i=dg_{0k}g_{0k}^{-1}=E_i\frac{dz}{z-s_i}.%
$$
So we obtained the holomorphic principal bundle ${\bf P}_\rho \to X$ on the
surface $X.$ The vector bundle associated to ${\bf P}_\rho \to X$, which we
denote by ${\bf E}_\rho \to X$ and call canonical, is not topologically
trivial. Its connexion is denoted by $\nabla$. The holomorphic sections of $%
{\bf E}_\rho $ are solutions of the equation 
$$
\nabla f=0\ \Longleftrightarrow \ df=\omega f.\eqno{(4.4)} 
$$

\ 

{\bf \ Theorem 4.1.}1) The system (4.4) has regular singularity at points $%
s_1,s_2,...,s_m$.

2) The Chern number $c_1($ ${\bf E}_\rho )$ of ${\bf E}_\rho \to X$ is equal
to 
$$
c_1({\bf E}_\rho )=\sum_{i=1}^mtr(E_i).\eqno{(4.5)} 
$$
The triple $\left( X,S,\rho \right)$ is called Riemann data, where $X$ is
a Riemann surface, $S\subset X$ denotes a finite subset of $X,$ $\rho :\pi
_1\left( X\setminus S,z_0\right) \to GL(n,{\bf C})$ is any representation
with trivial kernel.

\ 

{\bf Riemann-Hilbert monodromy problem for Riemann surfaces. }Let us find a
system of ODE 
$$
df=\omega f,\eqno{(4.6)} 
$$
on a Riemann surface $X$ for the given Riemann data $\left( X,S,\rho \right) 
$, where $S$ is the set of regular singular points of the system (4.6) and
its monodromy representation coincides with $\rho $.

\ 

{\bf Theorem 4.2 [R\"oh]}. For every Riemann data there exists a solution of
the Riemann-Hilbert problem for ODE's with regular singularity.

Let (4.6) be the regular system of ODE's which is induced by the
representation (4.1). By theorem 2.2 the fundamental matrix of solutions in
a neighbourhood of $s_j$ is
$$
\Phi _j\left( \tilde z\right) =U_j(z)(z-s_j)^{\Psi _j}(\tilde z-s_j)^{E_j}.
\eqno{(4.7)} 
$$
Here $\Psi _j$ are exponents of the solution space $\Re$ of the system
(4.6) and $E_j=\frac 1{2\pi i}\ln G_j,$ with eigenvalues $\mu _j^1,\mu
_j^2,...,\mu _j^n$ satisfying the conditons $0\leq Re\mu _j^i<1.$ The
numbers $\beta _j^i=\varphi_j^i+\mu _j^i$ will be called exponents of the
solution space $\Re $ at the point $s_j$ (or $j$-exponents).

{\bf Proposition 4.1 [Le], [Bl1]. }The system (4.6) is Fuchsian at $s_j$ if and
only if $\det U_j(s_j)\neq 0.$

{\bf Proposition 4.2.} If the system (4.5) is Fuchsian in a neighbourhood
of $s_j,$ then
$$
\omega _j=\frac{A_j}{z-s_j}dz, 
$$
where $A_i$ is a constant matrix with eigenvalues $\beta _j^i,i=1,...,n.$

{\bf Proof.} Indeed, suppose $\omega =A(z)dz,$ then
$$
A_j=\lim _{z\rightarrow s_j}\left( (z-s_j)A(z)\right) =\lim _{z\rightarrow
s_j}\left( (z-s_j)\frac{d\Phi _j(z)}{dz}\Phi ^{-1}(z)\right) = 
$$
$$
\lim _{z\rightarrow s_j}((z-s_j)\frac{dU_j(z)}{dz}U_j(z)+ 
$$
$$
+U_j(z)\left( \Psi _j+(z-s_j)^{\Psi _j}E_j(z-s_j)^{-\Psi _j}\right)
U_j^{-1}(z)= 
$$
$$
=U_j(s_j)(\Psi _j+E_j)U_j^{-1}(s_j). 
$$

Here $E_j=\lim _{z\rightarrow s_j}E_j(z),$ $E_j(z)=(z-s_i)^{\Psi
_j}E_j(z-s_i)^{-\Psi _j}.$ Therefore we obtain that $\psi _j^i+\mu _j^i$ are
eigenvalues of the matrix $A_j.$

Let the system (4.5) be Fuchsian. Transform the monodromy matrices $%
G_j,j=1,2,...m$ to upper-triangular form by some matrices $C_j.$ Assume that $%
\Psi _j,$ $j=1,2,...,m$ are diagonal integer valued matrices whose entries $%
\varphi _j^i$ satisfy the inequalities
$$
\varphi _j^1\geq \varphi _j^2\geq ...\geq \varphi _j^n. 
$$

Consider the local section $U_j(z)$ of the principal bundle ${\bf P}_\rho
\rightarrow X$ over $V_j\backslash s_j$ such that the corresponding $\Phi
_j(z)$ has the form (4.7).

The following proposition holds:

{\bf Proposition 4.3.} Every extension of ${\bf P}_\rho ^{\prime }\rightarrow
X_m$ to the points $s_j$ which is induced by a connexion $\nabla$ with at
most logarithmic singularities at $s_j,$ is determined by matrices $C_j$ and
$\Psi_j$ such that

1) $C_j^{-1}G_jC_j$ is upper triangular,

2) $\Psi _j=diag(\varphi _j^1,\varphi _j^2,...,\varphi _j^n),\varphi _j^i\in 
{\bf Z,}\varphi _j^1\geq \varphi _j^2\geq ...\geq \varphi _j^n.$

Extend ${\bf P}_\rho ^{\prime }\rightarrow X_m$ in a similar way to all
singular points. Denote by $C$ the collection $(C_1,C_2,...,C_m)$ and by
$\Psi$ the collection ($\Psi ^1,\Psi ^2,...,\Psi ^m),$ where
$\Psi ^j=(\varphi_j^1,\varphi _j^2,...,\varphi _j^n).$ Denote by
 ${\bf P}_\rho ^{C,\Psi }\rightarrow X$ the correspoding extension of the
bundle ${\bf P}_\rho^{\prime }\rightarrow X_m.$

The collection $C,\Psi$ is said to be admissible, if $C_j,\Psi ^j$ satisfy
1), 2) for every $j$.

{\bf Proposition 4.4. }There is a one-to-one correspodence between the set
of all Fuchsian systems of ODE's on the Riemann surface with prescribed
monodromy and the set $\left\{ H^0(X,O({\bf P}_\rho ^{C,\Psi })\right\} $ of
holomorphic sections of all admisible extensions of the principal bundle $%
{\bf P}_\rho ^{\prime }\rightarrow X_m.$

The proof of the last two propositions in the case when $X$ is the Riemann
sphere, is contained in [Bl3].

\section{Criterion of stability.}

Let $E\to X$ be a holomorphic vector bundle on a Riemann surface $X$, with $%
\deg E=k$ and rank$E=n$. The normalized Chern class of the vector bundle $E$
is defined by $\mu\left(E\right) =\frac kn.$

A bundle $E$ is called stable (resp. semi-stable) if for every proper
subbundle $F\subset E$, we have 
$$
\mu \left( F\right) <\mu \left( E\right),%
$$
(resp. 
$$
\mu \left( F\right) \leq \mu \left( E\right) ).%
$$

\ 

\ {\bf \ Properties of stable bundles:}

1) If $E\to X$ is semi-stable and $\left( n,k\right) =1,$ then $E$ is stable.

2) A line bundle $L\to X$ is stable.

3) Let $L\to X$ be a line bundle. $E\to X$ is stable if and only if $%
E\otimes L\to X$ is stable.

4) If $E\to X$ is stable, then $E$ is indecomposable and by the Riemann-Roch
theorem 
$$
\dim H^1\left( X,{\bf O}\left({\rm End}E\right) \right) =n^2(g-1)+1.%
\eqno{(5.1)} 
$$

The definition of stability has been given by D.~Mumford {\bf [Mum]} oriented
towards Riemann surfaces of negative curvature. For example, in case $g=0,$
stable bundle (in sense of Mumford) may be only line bundle $E\left(
k\right)$ (semi-stable are $E\left( k\right) ^{\oplus r})$, but there exist
stable bundles on the Riemann sphere, of rank more than one.

A criterion of stability for vector bundles on the Riemann sphere {\bf C}%
P$^1$ given by B.~Bojarski, will be reproduced below.

Let $\Im (n,k)$ be the space of all vector bundles on ${\bf CP}^1$ of rank $n
$ with Chern class $k$. By Grothendieck's theorem every bundle
$E\to{\bf CP}^1$ splits into direct sum (1.12) of line bundles.

The splitting type $K=(k_1,k_2,...k_n)\in {\bf Z}^n$ completely defines the
holomorphic structure of $E.$

{\bf Theorem 5.1. [Boj1] [Boj2].} 1) A vector bundle $E\to ${\bf C}P$^1$ is
stable if and only if $k_1-k_n\leq 1.$

2) The space of stable bundles is an open and dense subspace of $\Im(n,k)$.

If $E\to X$ is a holomorphic vector bundle over a Riemann surface of genus $%
g\geq 2$, then it does not split into the sum of line bundles but some
analogy exists {\bf [Gi]}.

A criterion of stability belongs to A.~Weil. In particular the following
theorem is true.

{\bf Theorem 5.2. [We]}. A topologically trivial vector bundle is stable if
and only if it corresponds to an irreducible unitary representation of the
fundamental group $\pi _1\left( X\right)$ 
$$
\rho :\pi _1\left( X\right) \to U(n).\eqno{(5.2)} 
$$

A generalization of Weil's theorem is Narasimhan and Seshadri theorem, which
is a criterion of stability for topologically nontrivial holomorphic vector
bundles.

{\bf Theorem 5.3.[N-S]. } A holomorphic vector bundle $E\to X$ of rank $n$
and with Chern class $k$ is stable if and only if it is induced from an
irreducible representation of the Fuchsian group, $\rho :\Gamma \to U(n),$
where $\Gamma $ is a group with $2g+1$ generators $%
a_1,b_1,a_2,b_2,...,a_g,b_g,c$ (where $g$ is genus of $X$), satisfying the
relations: 
$$
\prod_{i=1}^q\left[ a_i,b_i\right] c=1,\eqno{(5.3)} 
$$
$$
c^k=1,\eqno{(5.4)} 
$$
whereas the irreducible representation $\rho $ is 
$$
\rho \left( c\right) =\exp \left( -2\pi i\mu \left( E\right) \right) {\bf 1}%
. 
$$

Let $s\in X$ be a marked point and $k<0$ any integer. As it is known, there
exists a covering $\pi :{\bf H}\to X$, which is branched at the point $s$
and generators $\left( 5.3\right) $ of the uniformization group $\Gamma $
satisfying the conditions $\left( 5.4\right) -\left( 5.5\right) $ and ${\bf H%
}/\Gamma \cong X\setminus \left\{ s\right\} .$ So, $\Gamma $ contains a
finite cyclic subgroup with generator $c,$ it means that $\pi
^{-1}(S)\subset {\bf H}$ is the only fixed point of $\Gamma ,$ and $\Gamma $
is a central extension of $\pi _1(X):$ 
$$
1\to {\bf Z}_m\to \Gamma \to \pi _1(X)\to 1. 
$$

The group $\Gamma $ acts on the trivial bundle ${\bf H}\times {\bf C}^n\to X$
via 
$$
\left( z,c\right) \longmapsto \left( cz,\rho \left( c\right) v\right) .%
$$
This action gives us the holomorphic bundle ${\bf E}_\rho \to X\setminus \left\{
s\right\} .$

For the Riemann data $\left( X,\left\{ s\right\} ,\rho \right) $ by the
theorem 4.2 there exists the system of ODE's 
$$
df=\omega f,\eqno{(5.5)} 
$$
which has regular singular points and monodromy representation of the system
(5.5) coinciding with $\rho $. This means that $\omega $ is a connexion form
for the holomorphic bundle ${\bf E}_\rho \to X\setminus \left\{ s\right\} $.
The corresponding holomorphic connexion is denoted by $\nabla ^{\prime }.$

According to the construction of the 4-th section it is possible to extend $%
\left( {\bf E}_\rho ^{\prime },\nabla ^{\prime }\right) $ to ${\bf E}_\rho
\to X.$ By proposition 4.2 the connexion form $\omega $ in the neighbourhood 
of $s$ will be: 
$$
\omega =\mu \left( {\bf E}_\rho \right) {\bf 1}\frac{dz}{z-s}. 
$$

We obtain the following statement:

{\bf Theorem 5.4. }Let $E\rightarrow X$ be a stable holomorphic vector
bundle. Then there exists a Fuchsian connexion $\theta$ on $E$ which has
only one singular point.

The system (5.5) has apparent singular points, and their quantity will be
estimated if a local representation induced by $\rho $ is semi-simple.

A local representation induced by $\rho $ at a point $s_k\in S$ is defined
as follows: Let $U$ be a neighbourhood of $p$ that is biholomorphic to the
unit disk satisfying $U\cap S=\left\{ s_k\right\} .$ The injection $%
U\setminus \left\{ s_k\right\} \to M\setminus S$ induces a representation of 
$\pi _1\left( U\setminus \left\{ s_k\right\} \right) $ in $GL(n,{\bf C}).$
This is the local representation at $s_k\in S$ induced by $\rho$.

{\bf Theorem 5.5.[O]}. If the representation $\rho$ is irreducible and the
local representation at some point of $S$ induced by $\rho$ is semi-stable,
there exists a Fuchsian linear differential equation on $M$ which has at
most 
$$
1-n(1-g)+\frac{n(n-1)}2(m+2g-2)\eqno{(5.6)} 
$$
apparent singularities, where $m=$ card$S.$

This estimate will be obtained by calculating zeros of the Wronskian of our
system. Indeed,let $\Phi \left( z\right) $ be the fundamental system of
solutions of (5.5) and $f_1(z)$, $f_2(z)$, ... , $f_n(z)$ be any row. Denote
by $D$ the operator $\frac d{dz}.$ Consider the equation 
$$
\det \left( 
\begin{array}{cccc}
g(z) & f_1(z) & ... & f_n(z) \\ 
Dg(z) & Df_1(z) & ... & Df_n(z) \\ 
... & ... & ... & ... \\ 
D^ng(z) & D^nf_1(z) & ... & D^nf_n(z)
\end{array}
\right) =0,\eqno{(5.7)} 
$$
i.~e. 
$$
\omega _0(z)D^ng(z)+\omega _1(z)D^{n-1}g(z)+...+\omega _n(z)g(z)=0,%
\eqno{(5.8)} 
$$
where 
$$
\omega _0(z)=\det \left( 
\begin{array}{cccc}
f_1(z) & f_2(z) & ... & f_n(z) \\ 
Df_1(z) & Df_2(z) & ... & Df_n(z) \\ 
... & ... & ... & ... \\ 
D^nf_1(z) & D^nf_2(z) & ... & D^nf_n(z)
\end{array}
\right)  
$$
is the Wronskian of $f(z)=\left( f_1(z),f_2(z),...,f_n(z)\right) .$ We
denote it by $W(z).$ Write (5.8) as follows: 
$$
\frac 1{W(z)}(\omega _0(z)D^ng(z)+\omega _1(z)D^{n-1}g(z)+...+\omega
_n(z)g(z))=0.\eqno{(5.9)} 
$$

This equation is Fuchsian and its monodromy coincides with monodromy of the
system (5.5), but (5.9) has apparent singular points, they are zeros of the
Wronskian $W(z).$ It is clear, that the poles of $W(z)$ are singular points
of the system (5.5).

{\bf Collorary.} The connexion $\theta $ has at most
$$
n^2g-\frac{n(n-1)}2+1 
$$
apparent singular points.

Let $\left( X,S,\rho \right) $ be any Riemann data, where $S$ is empty and $%
\rho $ is an irreducible unitary representation, then (5.6) is an ODE with
apparent singular points$.$ The number $p$ of apparent singular points
will be estimated (by theorem 5.4): 
$$
p\leq n^2(g-1)+1. 
$$

Note that the right hand side of this inequality is the dimension of the
moduli space of stable holomorphic $n$-rank vector bundles on a Riemann
surface of genus $g.$

Every holomorphic bundle ${\bf E}$ has a canonical filtration {\bf [At-Bo]} 
$$
0=E_0\subset E_1\subset E_2\subset ...\subset E_r=E, 
$$
with $F_i=E_i/E_{i-1}$ semi-stable and 
$$
\mu (F_1)>\mu (F_2)>...>\mu (F_r). 
$$
$E$ is a direct sum of semi-stable bundles $F_1,F_2,...,F_p.$ Hence
$$
E\cong F_1\oplus F_2\oplus ...\oplus F_p.\eqno{(5.10)} 
$$

On the other hand, for any semi-stable vector bundle $V$ there exists a
filtration
$$
V=V_q\supset V_{q-1}\supset ...\supset V_1=\emptyset , 
$$
such that $V_1$ and $W_i=V_{i-1}/V_i$ are stable. Since $V_1$ is stable, it
follows that the cocycle $z^{-\mu (V_1)}$ defines $V_1.$ Analogically the
cocycle for $V_2$ is

$$
a_{V_2}=\left( 
\begin{array}{c}
(z-x_\infty )^{-\mu (W_1)}............* \\ 
0.............(z-x_\infty )^{-\mu (W_2)} 
\end{array}
\right) 
$$
and so on. Finally, for $V$ we have 
$$
a_V=\left( 
\begin{array}{c}
(z-x_\infty )^{-\mu (W_1)}............* \\ 
..... \\ 
0.............(z-x_\infty )^{-\mu (W_q)} 
\end{array}
\right). 
$$

By (5.10) it follows that 
$$
\Psi _E=diag(A_{F_1},...,A_{F_p}).\eqno{(5.11)} 
$$

This gives the proof of the following result.

{\bf Theorem 5.6. }Let $E\rightarrow X$ be a holomorphic vector bundle. Then
the cocycle (5.11) defines the given bundle up to an isomorphism.

The connexion $\omega$ agrees with the holomorphic structure on ${\bf E}_\rho
\to X$ and if $\omega _1$ is gauge equivalent to $\omega $, the bundle ${\bf %
E}_{\rho _1}\to X$ is holomorphically equivalent to ${\bf E}_\rho \to X.$
This gives possibility to describe holomorphic structures on the $C^\infty $%
-bundle ${\bf E}\to X.$

We consider therefore a fixed $C^\infty$ complex vector bundle ${\bf E}\to X$
of rank $n$ and Chern class $k$ and we denote the space of all holomorphic
structures on ${\bf E}$ by $\aleph (n,k)$.

Let Aut$(E)$ denote the group of automorphisms of $E$ which means that any
element of this group is locally a $C^\infty $ map of $X$ into $GL(n,{\bf C}%
).$ Then Aut$(E)$ acts on $\aleph (n,k)$ and the orbits by definition are
the isomorphism classes of complex analytic bundles on $X$ with Chern class
$k$ and rank $n.$

Suppose $F_i$ has rank $n_i$ and Chern class $k_i$ so that $n=\sum n_i$ and $%
k=\sum k_i.$ So, we have the sequence of rational numbers 
$$
\mu =\left( \frac{k_1}{n_1},...\frac{k_1}{n_1},\frac{k_2}{n_2},...,\frac{k_2%
}{n_2},...,\frac{k_r}{n_r},...\frac{k_r}{n_r}\right) , 
$$
which we call the type of ${\bf E.}$

All the holomorphic bundles of given type $\mu $ define a subspace $\aleph
_\mu $ of $\aleph (n,k).$ In particular, if all components of $\mu $ are
equal (hence are all equal to $\frac kn$), then $\aleph _\mu \cong \aleph
_{ss}$ is the semi-stable part of $\aleph (n,k).$

The codimension of $\aleph _\mu $ in $\aleph \left( n,k\right)$ is equal to 
$$
\sum_{\mu _i>\mu _j}^{}\left(\mu_i-\mu_j+\left(g-1\right)\right).
$$

If $g=0,$ we obtain the well known formula (theorem 1.3).

\section{Holomorphic bundles on CP$^1$}

Let $S=\left\{ s_1,s_2,...s_m\right\}$ be a set of marked points on ${\bf CP}%
^1$ and 
$$
df=\omega f\eqno{(6.1)}%
$$
be the system of ODE's which is induced by the representation 
$$
\rho :\pi _1\left( {\bf CP}^1\backslash S,z_0\right) \to GL(n,{\bf C}).%
\eqno{(6.2)}%
$$

By proposition 1.2 the fundamental matrix of solutions in the neighbourhood
of $s_j$ is 
$$
\Phi _j\left( \tilde z\right) =U_j(z)(z-s_j)^{\Psi _j}(\tilde z-s_j)^{E_j^{}}%
\eqno{(6.3)}, 
$$
where $\Psi _j$ are exponents of the solution space $\Re $ of the system
(6.1) and $E_j=\frac 1{2\pi i}\ln G_j.$ Again denote by $\mu _j^1,\mu
_j^2,...,\mu _j^n$ the eigenvalues of $E_j.$ Then eigenvalues $\beta
_j^i=\varphi _j^i+\mu _j^i$ of the matrix $\Psi _j+F_j$ are the exponents of
the solution space $\Re $ at the point $s_j$ (or $j$-exponents).

Let every singular point satisfy the condition $\det U_j(s_j)\ne 0.$ Then
the system (6.1) will be Fuchsian and therefore the 1-form $\omega $ will
have single poles. Denote 
$$
A_i=Res_{s_i}\omega , 
$$
then the system (6.1) will be written in the following form 
$$
df=\sum_{i=1}^m\frac{A_i}{z-s_i}dz\eqno{(6.4)}. 
$$
where the matrices $A_i$ satisfy the condition $\sum_{i=1}^mA_i=0.$

{\bf Proposition 6.1. [Le], [Bl1]} 1) The number 
$$
\beta =\sum_{j=1}^m\sum_{i=1}^n\beta _j^i 
$$
is integer and is at most 0.

2) The system (6.1) is Fuchsian if and only if $\beta=0.$

Our following discussion concerns Fuchsian systems. We emphasize this,
because after the Plemelj's work ``Problems in the sense of Riemann and
Klein'' it was believed that Hilbert's 21-st problem is solved for Fuchsian
systems, but A.~Bolibruch in 1989 gave an example of a representation of the
fundamental group $\pi _1({\bf CP}^1\backslash S,z_0)$ which does not give a
Fuchsian system.

{\bf Theorem 6.1.[Bl1]}. Let (${\bf CP}^1,S,\rho )$ be any Riemann data.

1) Case n=2. Then there exists a Fuchsian system where monodromy
representation coincides with $\rho$.

2) Case n=3. a) if $m=$Card$S=3$ or b) $m$ is arbitrary and $\rho $ is
irreducible then Hilbert's 21-st problem can be solved.

3) Case $n>3.$ For every $S,$ where Card$S>3,$ there exists a representation 
$\rho$ which does not induce a Fuchsian system.

Let us make few comments on this theorem.

1) Let n=2. There exists a representation of the fundamental group, which is
not monodromy representation for any Fuchsian system, even in case, when n=3.

\ {\bf Examples.} Let $s_1,s_2,s_3\in ${\bf C}P$^1$. Correspoding monodromy
matrices are 
$$
G_1=\left( 
\begin{array}{cc}
1 & c_1 \\ 
0 & 1 
\end{array}
\right) ,G_2=\left( 
\begin{array}{cc}
1 & c_2 \\ 
0 & 1 
\end{array}
\right) ,G_3=\left( 
\begin{array}{cc}
1 & -c_1-c_2 \\ 
0 & 1 
\end{array}
\right) , 
$$
where the numbers $c_1$, $c_2$ satisfy the condition $c_1c_2(c_1+c_2)\ne 0.$
There does not exist Fuchsian equation, whose monodromy group is generated
by these $G_1,G_2,G_3$ matrices. There always exists a Fuchsian system for
this representation. This follows from Dekkers' {\bf [D]} theorem.

2) For the system 
$$
df=\left( \left( 
\begin{array}{ccc}
0 & 1 & 0 \\ 
0 & z & 0 \\ 
0 & 0 & -z 
\end{array}
\right) \frac{dz}{z^2}+\frac 16\left( 
\begin{array}{ccc}
0 & 6 & 0 \\ 
0 & -1 & 1 \\ 
0 & -1 & 1 
\end{array}
\right) \frac{dz}{z+1}\right.%
$$
$$
\left.+\frac 12\left( 
\begin{array}{ccc}
0 & 0 & 2 \\ 
0 & -1 & -1 \\ 
0 & 0 & 1 
\end{array}
\right) \frac{dz}{z-1}+\frac 13\left( 
\begin{array}{ccc}
0 & -3 & -3 \\ 
0 & -1 & 1 \\ 
0 & -1 & 1 
\end{array}
\right) \frac{dz}{z-\frac 12}\right)f%
$$
$0,-1,1,\frac 12$ are regular singular points (The points $-1,1,\frac 12$
are Fuchsian singularities. However $0$ is not Fuchsian, as it is a pole of
order 2). Its monodromy representation $\rho$ is reduced. Therefore the
Fuchsian system does not exists, with these singular points and
representation $\rho.$

From theorem 6.1 n. 3) follows, in general case, that on the marked Riemann
sphere always exists a vector bunble, which has no Fuchsian connexions. So,
solvability of Hilbert's 21st problem depends on the conformal structure on
the Riemann sphere. The problem has topological character only for rank two
vector bundles. This case we consider below.

For every Fuchsian equation there exists a Fuchsian system (6.4), which has
the same singular point and monodromy. In particular, for the hypergeometric
equation 
$$
y''+\frac{\gamma +(\alpha +\beta +1)}{z(z-1)}y'-\frac{\alpha \beta }{%
z(z-1)}=0 
$$
the aforementioned system will be: 
$$
df=\left( \left( 
\begin{array}{cc}
0 & 0 \\ 
-\alpha \beta & -\gamma 
\end{array}
\right) \frac d{dz}+\left( 
\begin{array}{cc}
0 & 1 \\ 
0 & \gamma -(\alpha +\beta ) 
\end{array}
\right) \frac{dz}{z-1}\right) f 
$$

As we mentioned there exists a representation whose corresponding Fuchsian
equation does not exist without apparent singular points. But if
representation is irreducible, for the estimation of the quantity of
apparent singular points we may have more precise inequality than in theorem
5.5.

Consider the vector bundle $E^{C,\Psi }$ associated with the principal
bundle $P^{C,\Psi },$ which we introduced in section 4. Let $K^{C,\Psi
}=(k_1^{C,\Psi },k_2^{C,\Psi },...,k_n^{C,\Psi })$ be the splitting type of
the vector bundle $E^{C,\Psi }.$

The number
$$
\tau (E^{C,\Psi })=\sum_{i=1}^n(k_1^{C,\Psi }-k_i^{C,\Psi }) 
$$
is called the wieght of the bundle $\tau (E^{C,\Psi })$ and the number $%
\frac{\tau (E^{C,\Psi })}{rankE}$ is called the normalized weight.

{\bf Theorem 6.2 [Bl3]. }For any holomorphic vector bundle $E\rightarrow 
{\bf CP}^1$ with splitting type $K^{C,\Psi }=(k_1^{C,\Psi },k_2^{C,\Psi
},...,k_n^{C,\Psi })$ and for any points $S=\left\{ s_1,s_2,...s_m\right\} ,$
where
$$
m=\max (\max _{i=1,2,...,n-1}(k_i^{C,\Psi }-k_{i+1}^{C,\Psi }+2),3) 
$$
there exists a Fuchsian system with the singular points $s_1,s_2,...s_m$ and
with monodromy $\rho $ such that the following conditions hold:

1) $\rho $ is irreducible.

2) The canonical extension $E\rightarrow {\bf CP}^1$ of the vector bundle $%
E_\rho ^{\prime }\rightarrow {\bf CP}^1\backslash S,$ constructed by means
of the monodromy representation $\rho ,$ has the splitting type $k_1-\gamma
(E),$ $k_2-\gamma (E),...,k_n-\gamma (E),$ where $\gamma (E)$ is equal to
the integer part of the normalized Chern number of the vector bundle $E$.

3) For the vector bundle $E\rightarrow {\bf CP}^1$ there exists a
meromorphic connexion with at most logarithmic singularities and with an
irreducible monodromy whose number of singular points is equal $m$ and there
is no connexion with the properties mentioned above whose number of singular
points is less than $m$.

Let (${\bf CP}^1,S,\rho )$ be any Riemann data with card$S=m,$ $\rho :\pi
_1({\bf CP}^1\backslash S)\to GL(n,{\bf C)}$ be irreducible and let the type
of the canonical bundle induced from $\rho $ be $\left(
k_1,k_2,...k_n\right) .$ Denote by $l$ the quantity of the first equal
numbers, i.~e. $k_1=k_2=...=k_l.$ Under these conditions the quantity of
apparent singular points is at most 
$$
\frac{(m-2)n(n-1)}2-\sum_{i=1}^n(k_1-k_i)+1-l.\ \ \ {\bf [Bl2].} 
$$

Consequently, we obtain an estimate for the partial indices in the case
when $\rho $ is irreducible: 
$$
\sum_{i=1}^n(k_1-k_i)\leq \frac{(m-2)n(n-1)}2+1-l.\eqno{(6.5).} 
$$

If the equality is achieved it means that we have the ODE induced by $\rho $
which has no apparent singular points.

Consider the canonical extension $E\rightarrow {\bf CP}^1$ of the vector
bundle $E_\rho ^{\prime }\rightarrow {\bf CP}^1\backslash S$ induced by the
representation (6.2) of the Fuchsian system (6.4). The splitting type of $%
E\rightarrow {\bf CP}^1$ can be algorithmically calculated with the aid of
system (6.4) as follows (we are repeating here the nice argument from [Bl3]):

Consider the matrix $A_1$ of (6.4), with eigenvalues $\beta _1^1,$ $\beta
_1^2,...,$ $\beta _1^n.$ Let $\alpha _j^i=Re$ $\beta _j^i$ and $\alpha
_1^1\neq 0.$ Without loss of generality, we can assume that $\alpha _1^1\geq
1.$ (The case $\alpha _1^1\leq -1$ can be investigated in a similar way).

Consider the change
$$
g_1=T_1f\eqno{(6.6)} 
$$
of the dependent variable $f$ under the action of the constant nondegenerate
matrix $T_1$ such that the coefficient matrix $A_1^{\prime }=T_1A_1T_1^{-1}$
of the new system has Jordan normal form with the first eigenvalue equal
to $\beta _1^1.$ Recall that under the transformation (6.6) the system (6.1)
is transformed into the system
$$
\frac{dg_1}{dz}=A'(z)g_1\eqno{(6.7)} 
$$
where $A'=T_1AT_1^{-1}+\frac{dT_1}{dz}T_1^{-1}.$

Then consider the transformation
$$
g_2=(z-s_1)^Dg_1\eqno{(6.8)} 
$$
where $D=diag(-1,0,...,0).$ Under transformations (6.6), (6.8) our original
system (6.4) is transformed into system (6.1), which is Fuchsian at $%
s_1,s_2,...s_m$ with an additional apparent singular point $\infty .$%
Moreover, the eigenvalues of its coefficient matrix $A_1^{\prime \prime }$
at $s_1$ are equal to the following ones: $\beta _1^1-1,$ $\beta _1^2,...,$ $%
\beta _1^n.$

Using the procedure whose first step was described above, we can obtain
system (6.1), which is Fuchsian at $s_1,s_2,...s_m$ with the additional
apparent singularity at $\infty ,$ and whose coefficeient matrix
$\widetilde{A}_1$ has an eigenvalue ${\widetilde\beta}^1_1$ such that
$[Re\widetilde{\beta}_1^1]=0.$

In a similar way, we can obtain a system whose $\alpha _j^i$ at all points $%
s_1,s_2,...s_m$ are equal to zero. This system is Fuchsian at these points and
has one additional apparent singular point at $\infty .$ The transformation
matrix $T(z),$ which transforms our (6.4) system into this new one, is
meromorphic at $s_1,s_2,...s_m$ and $\infty ,$ holomorphically invertible off
these points, and can be calculated algorithmically on the basis of the
system (6.4).

Treat the matrix $T(z)$ as a transition function of some vector bundle on $%
{\bf CP}^1$ trivialized on ${\bf C}^1$ and a coordinate neighborhood $V_\infty$
of infinity. Then there exists a matrix $\Gamma (z)$ holomorphically invertible
in ${\bf C}^1$ and a matrix $U(z)$ holomorphically invertible in $V_\infty $
such that
$$
\Gamma (z)T(z)=z^KU(z), 
$$
where $K=diag(k_1,k_2,...,k_n),$ $k_1\geq k_2\geq ...\geq k_n.$

{\bf Proposition 6.2. }The collection of the numbers $-k_1,-k_2,...,-k_n$
coincides with the splitting type of the canonical extension $E\rightarrow 
{\bf CP}^1$ constructed by the monodromy of the Fuchsian system (6.4).

Consider the case n=2.

Let $i$-exponents for the solutions space $\Re $ be $\varphi _i^1,\varphi
_i^2 $ and let us assume $\varphi _i^1\geq \varphi _i^2$ (which is always
possible).

The number $\gamma _\omega =\sum_{i=1}^m\left( \varphi _i^1-\varphi
_i^2\right) $ is called Fuchsian weight for system (6.4).

Fix a representation $\rho$ for any Riemann data $({\bf C}P^1,S,\rho)$ and
denote by $\Omega_\rho$ the set of Fuchsian systems corresponding to this
data. The number $\gamma _\rho =\min _{\Omega _\rho }\gamma _\omega $ is
called Fuchsian weight for the representation $\rho .$

Let ${\bf E}\to{\bf C}P^1$ be the holomorphic vector bundle induced by the
representation $\rho $ and $(k_1,k_2)$ be its splitting type. Then 
$$
\gamma _\rho =k_1-k_2.%
$$

Every rank two holomorphic bundle on ${\bf C}P^1$ is holomorphically
equivalent to any bundle ${\bf F}\to{\bf C}P^1,$ which is obtained by the
extension of the bundle induced by an irreducible representation. So, for
every rank 2 holomorphic bundle there exists an irreducible connexion, which
is holomorphic except for the finite number $n_\omega $ of points, where it
has simple poles. Denote by $\Omega ^{\rm irr}$ the space of irreducible
Fuchsian connexions. Let 
$$
p=\min _{\omega \in \Omega ^{{\rm irr}}}n_\omega . 
$$

The identity 
$$
p=k_1-k_2+2\eqno{(6.9)}
$$
is satisfied.

From this follows the proposition.

\ 

{\bf Proposition 6.3. }A rank two vector bundle ${\bf F}\to ${\bf C}P$^1$ is
stable if and only if it is induced by the Gauss equation.

{\bf Proof.} Let the bundle be stable, then from (6.9) we obtain p=3. As we
mentioned, there exists an irreducible representation 
$$
\rho :\pi _1\left( {\bf CP}^1\backslash \left\{ s_1,s_2,s_3\right\}
,z_0\right) \to GL(n,{\bf C}), 
$$
from which ${\bf F}$ is induced. It means that there exists the Gauss
equation for every irreducible representation {\bf [Bl2]}. Therefore, ${\bf %
F}$ has Gauss connexion. Converse proposition is obtained from (6.9), taking
p=3.

\section{On the splitting type of a rank three vector bundle}

B.~Bojarski in [Boj3] has posed the question: whether the partial indices are
invariants of conformal transformation of the complex plane. Our following
reasoning is intended to move in this direction. 

Let $E\rightarrow {\bf CP}^1$ be a holomorphic vector bundle with
characteristic matrix $d_K,$ $K=(k_1,k_2,...,k_n),$ $k_1\geq k_2\geq ...\geq k_n.$
Consider the vector bundle End$E\rightarrow {\bf CP}^1.$ The characteriestic
matrix of this bundle is $d_K\otimes d_K^{-1}.$

Clearly
$$
\dim H^0\left( {\bf CP}^1,O(EndE)\right) =\sum_{k_i\geq k_j}(k_i-k_j+1), 
$$
$$
\dim H^1\left( {\bf CP}^1,O(EndE)\right)
=\sum_{k_i>k_j}(k_i-k_j-1)=\sum_{k_i>k_j}(k_i-k_j)+\frac 12n(n-1). 
$$
The number
$$
\nu (E)=\dim H^1\left( {\bf CP}^1,O(EndE)\right) -\frac 12n(n-1) 
$$
will be called reduced dimension of the deformation space of complex
structures of the bundle $E$.

For every holomorphic bundle $E\rightarrow {\bf CP}^1$ consider the
following diagram:
$$
\begin{array}{lllll}
k_1-k_2 & k_1-k_3 & ... & ... & k_1-k_n \\ 
k_2-k_3 & k_2-k_4 & ... & k_2-k_n &  \\ 
... & ... & ... & ... & ... \\ 
k_{n-1}-k_n &  &  &  &  
\end{array}
$$

The sum of the first row is the weight of the bundle
$E\rightarrow {\bf CP}^1,$ i.~e.
$$
\tau (E)=\sum_{i=1}^n(k_1-k_i)=nk_1-c_1(E)\eqno{(7.1)} 
$$
The sum of all elements of this diagram is equal to the reduced dimension of
the deformation space of the complex structures of the bundle $E$. From (7.1)
it follows:

1) The higher partial index of RHTP equals $k_1=\frac{\tau (E)}n+\frac{c_1(E)}%
n$ and therefore the sum of the normalized weight and normalized Chern
number of the vector bundle $E$ is an integer.

2) If $E$ is a rank two vector bundle, then $\tau (E)=\nu (E).$

{\bf Theorem 7.1. }Let $E\rightarrow {\bf CP}^1$ be a rank three vector
bundle. Then the splitting type of $E$ is given by
$$
k_1=\frac 13(c_1(E)+\tau (E)), 
$$
$$
k_2=\frac 13c_1(E)-\frac 23\tau (E))+\frac 12\nu (E), 
$$
$$
k_3=\frac 13(c_1(E)+\tau (E))-\frac 12\nu (E). 
$$

{\bf Theorem 7.2. }Let $E\rightarrow {\bf CP}^1$ be a rank two vector
bundle. Then
$$
k_1=\frac 12(c_1(E)+\nu (E)), 
$$
$$
k_2=\frac 12(c_1(E)-\nu (E)). 
$$

\bibliographystyle{alpha}

\ 

Institute of Cybernetics

Georgian Academy of Sciences

Sandro Euli st. 5

Tbilisi 380086

Republic of Georgia

email: giorgadze@@rmi.acnet.ge

\end{document}